\documentclass[11pt]{amsart}
\usepackage{amsmath,amssymb,color,cite,verbatim}
\newtheorem{lemma}[equation]{Lemma}
\newtheorem{cor}[equation]{Corollary}

\theoremstyle{definition}
\newtheorem*{rmk}{Remark}

\newcommand{\F}{\mathbb{F}}
\newcommand{\bP}{\mathbb{P}}

\DeclareMathOperator{\ord}{ord}

\usepackage[colorlinks,pdftex,pagebackref, bookmarks=false]{hyperref}

\begin{document}

\title{A note on the paper arXiv:2207.13335}

\author{Michael E. Zieve}
%\address{
%  Department of Mathematics,
%  University of Michigan,
%  530 Church Street,
%  Ann Arbor, MI 48109-1043 USA
%  }
%\email{zieve@umich.edu}
%\urladdr{http://www.math.lsa.umich.edu/$\sim$zieve/}

%\date{\today}

\begin{abstract}
We show that all of the ``new" permutation polynomials in the
recent paper arXiv:2207.13335 are in fact known.  We also present a new type of question in this area.
\end{abstract}

\maketitle

%#######################################################################
%#######################################################################

%\section{Introduction}

A polynomial $f(X)\in\F_q[X]$ is called a \emph{permutation polynomial} if the function $c\mapsto f(c)$ permutes $\F_q$.  There has been much recent interest in the class $\mathcal F_q$ of permutation polynomials over $\F_{q^2}$ of the form $X^r A(X^{q-1})$.
The main reason this particular form is special is that a procedure in \cite{ZR} shows how to produce permutation polynomials in $\mathcal F_q$ from any prescribed permutation polynomial (or more generally, permutation rational function) over $\F_q$.  A second special feature of $\mathcal F_q$ is that there is a known method to use any prescribed polynomial in $\mathcal F_q$ in order to produce arbitrarily many other polynomials in $\mathcal F_q$.

The recent paper \cite{SGZWL} purports to produce new classes of permutation polynomials.  Here we show that all the permutation polynomials in that paper are obtained by applying the above-mentioned method to some well-known permutation polynomials in $\mathcal F_q$.

In more detail, the following special case of \cite[Lemma~1.2]{Zlem} reduces the study of permutation polynomials over $\F_{q^2}$ to the study of permutations of the group $\mu_{q+1}$ of all $(q+1)$-th roots of unity in $\F_{q^2}$:

\begin{lemma}\label{old}
Write $f(X):=X^r A(X^{q-1})$ where $r$ is a positive integer, $q$ is a prime power, and $A(X)\in\F_{q^2}[X]$.  Then $f(X)$ permutes\/ $\F_{q^2}$ if and only if $\gcd(r,q-1)=1$ and $g(X):=X^r A(X)^{q-1}$ permutes $\mu_{q+1}$.
\end{lemma}

The paper \cite{ZR} shows how to produce permutations of $\mu_{q+1}$ of the form $X^r A(X)^{q-1}$ from any prescribed rational function in $\F_q(X)$ which permutes $\bP^1(\F_q):=\F_q\cup\{\infty\}$.  More information about this procedure is given in \cite{Zx}, and full details appear in the forthcoming paper \cite{pp-struc1}.  In the present paper we only use some known permutations of $\mu_{q+1}$, so we need not say more about this procedure here.

The method to produce new permutations of $\mu_{q+1}$ from a known permutation is encoded in the following variant of \cite[Cor.~1]{DeZh}, whose proof is identical to that of \cite[Cor.~1]{DeZh}:

\begin{lemma}\label{SR}
Let $q$ be a prime power, let $r$ be an integer, and let $s_1,\dots,s_m$ and $t_1,\dots,t_m$ be positive integers for some $m\ge 0$.  Pick $A(X)\in\F_{q^2}[X]$ and write $B(X):=A(X)\cdot\prod_{i=1}^m \sum_{j=0}^{s_i} X^{jt_i}$.  Then $X^r B(X)^{q-1}$ permutes $\mu_{q+1}$ if and only if $X^{r-\sum_{i=1}^m s_it_i} A(X)^{q-1}$ permutes $\mu_{q+1}$ and, for each $i$, we have $\gcd(s_i+1,q)=1$ and $(q+1)/\gcd(t_i,q+1)$ is coprime to $s_i+1$.
\end{lemma}

\begin{rmk}
Other results along the lines of Lemma~\ref{SR} are
\cite[Lemma~4.2]{ZHF} and \cite[Lemma~4.1]{XCP}.  A generalization of all of these results, based on \cite[Thm.~5.1]{ZR}, will appear in \cite{pp-struc1}.
\end{rmk}

The above results have the following immediate consequence.

\renewcommand{\theenumi}{\thethm.\arabic{enumi}}
\renewcommand{\labelenumi}{(\thethm.\arabic{enumi})}

\begin{cor}\label{main}
Let $q$ be a prime power, and pick an integer $v$ and a polynomial $D(X)\in\F_{q^2}[X]$ such that $X^v D(X)^{q-1}$ permutes $\mu_{q+1}$.
Let $m$ be a nonnegative integer, and let $s_1,\dots,s_m$ and $t_1,\dots,t_m$ be positive integers such that for each $i$ we have $\gcd(s_i+1,q)=1$ and $(q+1)/\gcd(t_i,q+1)$ is coprime to $s_i+1$.
\begin{enumerate}
\item\label{31} If $r$ is a positive integer such that $\gcd(r,q-1)=1$ and $r\equiv v+\sum_{i=1}^m s_it_i\pmod{q+1}$ then $X^r B(X^{q-1})$ permutes\/ $\F_{q^2}$ where $B(X):=D(X)\prod_{i=1}^m \sum_{j=0}^{s_i} X^{jt_i}$. 
\item\label{32} If $B(X)\in\F_{q^2}[X]$ satisfies $D(X)=B(X)\prod_{i=1}^m\sum_{j=0}^{s_i} X^{jt_i}$, and $r$ is a positive integer such that $\gcd(r,q-1)=1$ and $r\equiv v-\sum_{i=1}^m s_it_i\pmod{q+1}$, then $X^r B(X^{q-1})$ permutes\/ $\F_{q^2}$.
\end{enumerate}
\end{cor}

In order to give explicit examples of the above result we now exhibit certain permutations $X^v D(X)^{q-1}$ of $\mu_{q+1}$.  Here we write $\ord_2(n)$ for the largest nonnegative integer $i$ such that $2^i\mid n$.

\begin{lemma}\label{pp}
Let $k$ and $\ell$ be positive integers, and write $q:=2^k$ and $Q:=2^\ell$.  Then $X^{Q+1} D(X)^{q-1}$ permutes $\mu_{q+1}$ when either of the following hold:
\begin{enumerate}
\item\label{41} $\ord_2(\ell)\le\ord_2(k)$ and $D(X):=X^{Q+1}+X+1$;
\item\label{42} $\ord_2(\ell)\ne\ord_2(k)$ and $D(X):=X^Q+X+1$.
\end{enumerate}
\end{lemma}

\begin{rmk}
This is a very special case of \cite[Thm.~1.3]{DZquad}.
It is a reformulation of \cite[Thm.~1 and 2]{LHnewtri}.  Other special cases of Lemma~\ref{pp} in the literature (stated in different but equivalent ways) are \cite[Thm.~3.1]{WZZ},
\cite[Cor.~3.8 and 3.10--3.12]{WYDM}, and \cite[Thm.~1.1]{Zx}.
\end{rmk}

Combining the above result with the special case of Corollary~\ref{main} in which $m=1$ and $s_1=2$ yields the following:

\begin{cor}\label{main2}
Let $k,\ell,t$ be positive integers with $k$ even, and write $q:=2^k$ and $Q:=2^\ell$.  Pick $D(X)\in\F_2[X]$ such that either \eqref{41} or \eqref{42} holds.
\begin{enumerate}
\item\label{51} If $r$ is a positive integer such that $\gcd(r,q-1)=1$ and $r\equiv Q+1+2t\pmod{q+1}$ then $X^r B(X^{q-1})$ permutes\/ $\F_{q^2}$ where $B(X):=D(X)\cdot (X^{2t}+X^t+1)$.
\item\label{52} If either $\ell$ is even and \eqref{41} holds or $\ell$ is odd and \eqref{42} holds then $B(X):=D(X)/(X^2+X+1)$ is in\/ $\F_{q^2}[X]$, and if $r$ is a positive integer such that $\gcd(r,q-1)=1$ and $r\equiv Q-1\pmod{q+1}$ then $X^r B(X^{q-1})$ permutes\/ $\F_{q^2}$.
\end{enumerate}
\end{cor}

Theorems 1, 3, 5 and 6 of \cite{SGZWL} are special cases of Corollary~\ref{main2}. 
% Actually the theorems in \cite{SGZWL} are false as stated, since they involve "polynomials" with negative exponents.  In what I wrote here, I implicitly assume that they intended to restrict to polynomials that are actually polynomials.  Also my statements of the result differ from theirs in that I have removed superfluous hypotheses and I have stated conditions on ord_2(k) and ord_2(ell) rather than equivalent conditions on gcd(Q-1,q+1) and gcd(Q+1,q+1).
Theorems 2 and 4 of \cite{SGZWL} are obtained from special cases of Corollary~\ref{main2} by composing on the right with $X^e$ for some positive integer $e$ such that $\gcd(e,q^2-1)=1$, and then reducing the composition mod $X^{q^2}-X$.  Other instances of Corollary~\ref{main2} are items (4) and (5) in \cite[Example~1]{DeZh}.

\begin{comment}
q=2^k, k even
s=2

(Q+1,q+1)=1 when ord_2(k)\ne ord_2(ell)
(Q+1,q-1)=1 when ord_2(ell)\ge ord_2(k)
(Q-1,q+1)=1 when ord_2(k)\ge ord_2(ell)

\end{comment}

One can write down arbitrarily many classes of permutation polynomials over $\F_{q^2}$ by applying Corollary~\ref{main} to any prescribed class of permutations of $\mu_{q+1}$ of the form $X^r D(X)^{q-1}$, using any choices of $s_i$'s and $t_i$'s.  However, as demonstrated in \cite{SGZWL}, this usually yields permutation polynomials with complicated expressions, and it is not clear how these enhance our understanding of the topic.  Instead, it would be interesting to systematically study the algebraic forms of all the permutation polynomials in $\mathcal F_q$ produced by Corollary~\ref{main} (and its variants) from some prescribed class of permutations of $\mu_{q+1}$ of the form $X^r D(X)^{q-1}$, and examine which of these permutation polynomials have an unexpectedly nice algebraic form, such as having few terms.
This is an important new type of question in this area, and we encourage people working on the topic to consider instances of it.

As a final remark, we urge readers to be cautious when using the list of all known permutation polynomials over $\F_{2^{2k}}$ in \cite[Thm.~7]{SGZWL}, since that list omits most of the known examples, such as $X^n$, Dickson polynomials, linearized polynomials, and many others.

%#######################################################################
%#######################################################################
%#######################################################################

%\bibliographystyle{plain}

\end{document}